\documentclass[graybox]{svmult}

\usepackage{type1cm}        % activate if the above 3 fonts are
\usepackage{makeidx}         % allows index generation
\usepackage{graphicx}        % standard LaTeX graphics tool
\usepackage{multicol}        % used for the two-column index
\usepackage[bottom]{footmisc}% places footnotes at page bottom
\usepackage{todonotes}

\usepackage{newtxtext}       % 
\usepackage{newtxmath}       % selects Times Roman as basic font
\usepackage[caption=false]{subfig}
\usepackage{tikz}
\usetikzlibrary{decorations.pathreplacing}
\usepackage{placeins}
\newcommand{\diagentry}[1]{\text{\fboxsep.5ex\fbox{$#1$}}}
\newcommand{\jump}[1]{\left\llbracket #1 \right\rrbracket}

\newcommand{\newstablong}{\text{Domain-of-Dependence stabilization}}
\newcommand{\newstab}{\text{DoD stabilization} }
\newcommand{\newstababbrv}{\text{DoD} }
\makeindex             % used for the subject index
\begin{document}

\title*{Monotonicity considerations for stabilized DG cut cell schemes for the unsteady advection equation}
\titlerunning{Monotonicity considerations for stabilized DG cut cell schemes}
% Use \titlerunning{Short Title} for an abbreviated version of
% your contribution title if the original one is too long
\author{Florian Streitb\"urger, Christian Engwer, Sandra May, and Andreas N\"u\ss ing}

% Use \authorrunning{Short Title} for an abbreviated version of
% your contribution title if the original one is too long
\institute{Florian Streitb\"urger, \at TU Dortmund University, Vogelpothsweg 87, 44227 Dortmund, Germany, \email{florian.streitbuerger@math.tu-dortmund.de}
\and Christian Engwer, Andreas N\"u\ss ing, \at University of M\"unster, Einsteinstraße 62, 48149 M\"unster,  Germany, \email{christian.engwer@uni-muenster.de, andreas.nuessing@uni-muenster.de}
\and Sandra May, \at Technical University of Munich, Boltzmannstraße 3,
85748 Garching bei München, Germany, \email{sandra.may@ma.tum.de}}
%
% Use the package "url.sty" to avoid
% problems with special characters
% used in your e-mail or web address
%
\maketitle

\abstract{
For solving unsteady hyperbolic conservation laws on cut cell meshes, the so called \textit{small cell problem} is a big issue: one would like to use a time step that is chosen with respect to the background mesh and use the same time step on the potentially arbitrarily small cut cells as well. For explicit time stepping schemes this leads to instabilities.
In a recent preprint [arXiv:1906.05642], we propose penalty terms for stabilizing a DG space discretization to overcome this issue for the unsteady linear advection equation. The usage of the proposed stabilization terms results in stable schemes of
first and second order in one and two space dimensions. In one dimension, for piecewise constant data in space and explicit Euler
in time, the stabilized scheme can even be shown to be monotone. In this contribution, we will examine the conditions for monotonicity in more detail.
}
%\newpage
\section{A stabilized DG cut cell scheme for the unsteady advection equation}
We consider the time dependent linear advection problem on a cut cell mesh. 
In \cite{May2019}, we propose new stabilization terms for a cut cell discontinuous Galerkin (DG) discretization in two dimensions with piecewise linear polynomials. In the following we will refer to this as \textit{\newstablong{}}, abbreviated by DoD stabilization.

While the usage of finite element schemes on embedded boundary or cut cell meshes has become
increasingly popular for elliptic and parabolic problems in recent years, only very little work has been done for hyperbolic problems. The general challenge is that cut cells can have
various shapes and may in particular become arbitrarily small. 
Special schemes have been developed to guarantee stability. Perhaps the most prominent approach for elliptic and parabolic problems is the ghost penalty stabilization \cite{Burman},
%sometimes referred to as the cutFEM method, 
which regains coercivity, independent of the cut size.

For hyperbolic conservation laws the problems caused by cut cells are partially of different nature. One major challenge is that standard explicit schemes are not stable on the arbitrarily small cut
cells when the time step is chosen according to the cell size of the background mesh.
This is what is often called the small
cell problem. Adapting the time step size to the cut size is infeasible, as there is no lower bound on the cut size.
An additional complication is the fact that there is typically no concept of coercivity
that could serve as a guideline for constructing stabilization terms.

In \cite{May2019}, we consider the small cell problem for the unsteady linear advection equation. We propose a stabilization of the spatial discretization, which uses a standard DG scheme with upwind flux, that makes explicit time stepping stable again. Our penalty terms are designed to restore the correct domains of dependence of the cut cells and their outflow neighbors (therefore the name DoD stabilization), similar to the idea behind the $h$-box scheme \cite{hbox_2003} but realized in a DG setting using penalty terms. 
In one dimension, we can prove $L^1$-stability, monotonicity, and TVD (total variation diminishing) stability for the stabilized scheme of first order using explicit Euler in time. For the second-order scheme, we can show a TVDM (TVD in the means) result if a suitable limiter is used.

In this contribution, we will focus on the monotonicity properties in one dimension for the first-order scheme and examine them in more detail. In particular, we will show that a straight-forward adaption of the ghost-penalty approach \cite{Burman} to the unsteady transport problem, as proposed in \cite{Massing2018} for the steady problem, cannot ensure monotonicity. Further, we will examine the parameter that we use in our new \newstab{} in more detail than done in \cite{May2019}.
%\textcolor{blue}{
%Using a method of lines approach, we start with a standard upwind DG discretization for the
%background mesh and add penalty terms that stabilize the solution on small cut cells in a conservative
%way. Then, one can use explicit time stepping, even on cut cells, with a time step length that is
%appropriate for the background mesh. In one dimension, we show monotonicity of the proposed
%scheme for piecewise constant polynomials and total variation diminishing in the means stability
%for piecewise linear polynomials. We also present numerical results in one and two dimensions that
%support our theoretical findings.
%}
\section{Problem setup for piecewise constant polynomials}
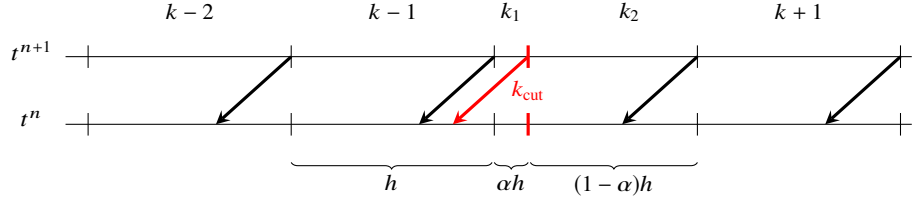
\begin{figure}
\begin{tikzpicture}[scale = 1.5, decoration = brace]
\draw (0,0.5) -- (7.5,0.5);
\draw (0,1.1) -- (7.5,1.1);
\draw (0.2,0.4) -- (0.2,0.6);
\draw (2,0.4) -- (2,0.6);
\draw (3.8,0.4) -- (3.8,0.6);
\draw [very thick, red] (4.1,0.4) -- (4.1,0.6);
\draw (5.6,0.4) -- (5.6,0.6);
\draw (7.4,0.4) -- (7.4,0.6);
%\draw (9.2,0.3) -- (9.2,0.7);
\draw (0.2,1.0) -- (0.2,1.2);
\draw (2,1.0) -- (2,1.2);
\draw (3.8,1.0) -- (3.8,1.2);
\draw [very thick, red] (4.1,1.0) -- (4.1,1.2);
\draw (5.6,1.0) -- (5.6,1.2);
\draw (7.4,1.0) -- (7.4,1.2);
%\draw (9.2,1.8) -- (9.2,2.2);

\draw [very thick,->,>=stealth] (2,1.1) -- (1.333,0.5);
\draw [very thick,->,>=stealth] (3.8,1.1) -- (3.133,0.5);
\draw [very thick,red,->,>=stealth] (4.1,1.1) -- (3.433,0.5);
\draw [very thick,->,>=stealth] (5.6,1.1) -- (4.933,0.5);
\draw [very thick,->,>=stealth] (7.4,1.1) -- (6.733,0.5);

%\draw () -- ();
\node at (1.1,1.5) {\small $k-2$};
\node at (2.9,1.5) {\small $k-1$};
\node at (3.95,1.5) {\small $k_1$};
\node at (5.0,1.5) {\small $k_2$};
\node at (6.5,1.5) {\small $k+1$};
%\node at (8.3,2.5) {\small $k+2$};
\node at (4.1,0.8) {\textcolor{red}{$k_\text{cut}$}};
\draw[decorate, yshift=-4ex] (4.1,0.7) -- node[below=0.4ex] {$\alpha h$} (3.8,0.7);
\draw[decorate, yshift=-4ex] (5.6,0.7) -- node[below=0.4ex] {$(1-\alpha) h$} (4.12,0.7);
\draw[decorate, yshift=-4ex] (3.78,0.7) -- node[below=0.4ex] {$ h$} (2,0.7);

%\draw [red] (2.9,0) -- (4.4,0);
%\draw [red] (2.9,-0.2) -- (2.9,0.2);
%\draw [red] (4.4,-0.2) -- (4.4,0.2);

\node at (-0.3,0.55) {$t^{n}$};
\node at (-0.3,1.15) {$t^{n+1}$};

\end{tikzpicture}
\caption{Domains of dependence for the solution at time $t^{n+1}$ for the considered model problem for a time step with length $\Delta t= \frac{\lambda}{\beta}h$ with $\lambda=\frac{1}{3}$ and $\beta =1$. }\label{fig:MP}
\end{figure}

For the purpose of a theoretical analysis with focus on monotonicity, we will consider piecewise constant polynomials in 1D. We use the interval $I=[0,1]$ and assume the velocity $\beta > 0$ to be constant. The time dependent linear advection equation reads as
\begin{equation}\label{eq: 1d adv}
u_t(x,t) + \beta u_x(x,t) = 0 \text{ in } I \times (0,T),
\end{equation}
with initial data $u(x,0) =u_0(x)$ and periodic boundary conditions.
We discretize the interval $I$ in N equidistant cells $I_j=[x_{j-\frac{1}{2}},x_{j+\frac{1}{2}}]$ with cell length $h$. Then, we split one cell, the cell $k$, into a pair of two cut cells using  the volume fraction $\alpha \in (0,\frac{1}{2}]$, see figure \ref{fig:MP}: The first cut cell, which we call $k_1$, has length $\alpha h$, the second cut cell, which we call $k_2$, has length $(1-\alpha)h$. Therefore, cell $k_1$ corresponds to a small cut cell, which induces instabilities, if $\alpha \ll \frac{1}{2}$.

We define the function space
\begin{equation}\label{eq:Funktionenraum}
V_h^0(I):=\left\{v_h \in L^2(I)\  \vline \ v_{h}|_{I_j} \in \mathcal{P}^0, \: j = 1,\ldots,N \right\},
\end{equation}
with $\mathcal{P}^0$ being the function space of constant polynomials. The semidiscrete scheme, which uses the standard DG scheme with an upwind flux in space and is not yet discretized in time, is given by: Find $u_h \in V_h^0(I)$ such that 
\begin{equation}\label{eq:semidiskret}
\int_I d_t u_h(t) \: w_h \: dx + a^{\text{upw}}_h(u_h(t),w_h)=0,\quad \forall w_h \in V_h^0(I),
\end{equation}
with the bilinear form defined as 
\begin{equation*}
a^\text{upw}_h(u_h,w_h)=\sum_{j=1}^N \beta u_h(x_{j+\frac{1}{2}}^-)\jump{w_h}_{j+\frac{1}{2}},
\end{equation*}
and the jump being given by
\[
\jump{w_h}_{j+\frac{1}{2}} = w_h(x_{j+\frac{1}{2}}^-)-w_h(x_{j+\frac{1}{2}}^+), \quad x_{j+\frac{1}{2}}^\pm = \lim_{\varepsilon\rightarrow 0^+}x_{j+\frac{1}{2}}\pm \varepsilon.
\]

We use explicit Euler in time. Then, \eqref{eq:semidiskret} results in the global system
\begin{equation}\label{eq:globalsystem}
\mathcal{M}u^{n+1} = \mathcal{B}u^n.
\end{equation}
Here, $u^n = [u_1^n,\ldots,u_N^n]^T$ is the vector of the piecewise constant solution at time
$t^n$ and $\mathcal{M}$ is the global mass matrix. Note that $\mathcal{M}$ is a diagonal matrix with positive diagonal entries. Finally, the global system matrix $\mathcal{B}$ is given by $\mathcal{B} = \mathcal{M} - \Delta t \mathcal{A}$ with $\mathcal{A}u^n$ corresponding to the discretization of the bilinear form
$a^\text{upw}_h$ at time $t^n$.

For a standard equidistant mesh, the scheme \eqref{eq:globalsystem} is stable for
$0 < \lambda < 1$ with the CFL number $\lambda$ being given by
\begin{equation}\label{eq:cfl}
\lambda = \frac{\beta \Delta t}{h}.
\end{equation}
Our goal is to make the scheme stable for the mesh containing the cut cell pair for
$0 < \lambda < \frac{1}{2}$, independent of $\alpha$. The reduced CFL condition is due to the fact that we will only stabilize cut cell $k_1$, and not the bigger cut cell $k_2$.

To illustrate one interpretation of the small cell problem that we need to overcome, we refer to figure \ref{fig:MP}. There, we determine the exact solution at time $t^{n+1}$, based on piecewise constant data at time $t^n$, by tracing back characteristics. For standard cells $j$, such as $k-1$, the domain of dependence of $u^{n+1}_j$ only includes cells $j$ and $j-1$. For the outflow neighbor of the small cut cell $k_1$, the cell $k_2$, however,
$u^{n+1}_{k_2}$ depends on $u^n_{k-1}$, $u^n_{k_1}$, and $u^{n}_{k_2}$. 
The issue is that standard DG schemes such as \eqref{eq:semidiskret} only provide information from direct neighbors. We will see that the proposed stabilization that ensures monotonicity will also fix this problem. We will return to this specific interpretation of the small cell problem in section \ref{sec:4}, when discussing the proper choice of the penalty parameter in the stabilization.
\section{Monotonicity considerations for two different stabilization terms}\label{sec:3}
In the following, we will examine the monotonicity properties of
different stabilizations. A monotone scheme guarantees in particular that
$\min_j u_j^0 \le u^n \le \max_j u_j^0$ for all times $t^n$. We will use the following definition of a monotone scheme.
\begin{definition}\label{Monoton}
A method
$ u^{n+1}_j = H(u^n_{j-i_L},u^n_{j-i_L+1},...,u^n_{j+i_R}) $
is called {\em monotone}, if for all $j$ there holds for every $l$ with $-i_L\le l\le i_R$
\begin{equation}\label{Def_monotone_coeff}
\frac{\partial H}{\partial u_{j+l}}(u_{j-i_L},...,u_{j+i_R})\geq 0.
\end{equation}
\end{definition}
For the linear scheme \eqref{eq:globalsystem} this implies that all coefficients of $\mathcal{B}$ need to be non-negative. This is due to the fact that $\mathcal{M}$ is a non-negative diagonal matrix. On an equidistant mesh, the scheme \eqref{eq:globalsystem} is monotone for $0 < \lambda < 1$.

We will compare the entries of the matrix $\mathcal{B}$ for three different cases: 
the unstabilized case, a stabilization in the spirit of the ghost-penalty method \cite{Burman}, and the \newstab \cite{May2019} that we propose.
\subsubsection*{Unstabilized case}
In this case, the matrix $\mathcal{B}$ is given by
\begin{equation*}
  \mathcal{B}=
    \left(
    \begin{smallmatrix}
      h-\tau & 0 & \cdots&  & & \cdots & 0 & \tau\\
      \tau & h-\tau & 0 & & & & & 0\\
      0 & \ddots & \ddots & & & & &\vdots\\
      \vdots & & \tau& \diagentry{\alpha h-\tau} & & & &\\
      & & & \tau& \diagentry{(1-\alpha)h-\tau} & & &\vdots \\
      \vdots & & & & \ddots& \ddots& & 0\\
      &&&&&&&\\
      0 &\cdots& & &\cdots&0&\tau& h-\tau\\
    \end{smallmatrix}\right),
\end{equation*}
with $\tau:=\beta\Delta t > 0$. We therefore focus on the diagonal entries.
On standard cells, and on cell $k_2$, the entries $h-\tau$ and $(1-\alpha)h-\tau$ are non-negative due to the CFL condition $\beta \Delta t = \lambda h$ if the reduced CFL condition $0 < \lambda \le \frac{1}{2}$ is used. 
On the small cut cell $k_1$, the entry $\alpha h - \tau$ is clearly negative for $\alpha < \lambda$, which is the case that we are interested in.
\subsubsection*{Ghost-penalty stabilization}
We first consider the option of using the ghost penalty method for stabilization, 
an approach that is, e.g., used in \cite{Massing2018} for stabilizing the steady advection equation. 
Adapting the stabilization to our model mesh (compare figure \ref{fig:MP}) changes the formulation of \eqref{eq:semidiskret} to: Find $u_h \in V_h^0(I)$ such that 
\begin{equation}\label{eq:gp_formulation}
\int_I d_t u_h(t) \: w_h \: dx  + a^{\text{upw}}_h(u_h(t),w_h)+J_h^{GP}(u_h,w_h)=0,\quad \forall w_h \in V_h^0(I),
\end{equation}
with 
\begin{equation}\label{eq:gpstab}
J_h^{GP}=\beta\eta_1\jump{u_h}_{k-\frac{1}{2}}\jump{w_h}_{k-\frac{1}{2}}+\beta\eta_2\jump{u_h}_{k_\text{cut}}\jump{w_h}_{k_\text{cut}}.
\end{equation}
As a result, the matrix $\mathcal{B}$ in \eqref{eq:globalsystem} is modified in the following way
\begin{equation*}
  \mathcal{B}_{GP}=
    \left(
    \begin{smallmatrix}
      h-\tau & 0 & \cdots& & & & & \cdots & 0 & \tau\\
      \tau & h-\tau & 0 & & & & & & & 0\\
      0 & \ddots & \ddots & \ddots & & & & & &\vdots\\
      &&&&&&&&\\
      &&&&&&&&\\
      & & & \tau & \diagentry{h-\tau(1-\eta_1)} & -\tau \eta_1 & & &\\
      \vdots & && 0 & \tau(1-\eta_1)& \diagentry{\alpha h-\tau+\tau\eta_1+\tau\eta_2} & -\tau\eta_2& & &\\
      & & & & 0 & \tau(1-\eta_2)& \diagentry{(1-\alpha)h-\tau(1-\eta_2)} & & &\vdots \\
      &&&&&&&&&\\
      &&&&&&&&&\\
      \vdots & & & & & & \ddots& \ddots& & 0\\
      &&&&&&&&&\\
      0 &\cdots& & & & &\cdots&0&\tau& h-\tau\\
    \end{smallmatrix}\right).
\end{equation*}
Our goal is to determine the parameters $\eta_1$ and $\eta_2$ such that every entry of $\mathcal{B}_{GP}$ is non-negative. The two entries on the first superdiagonal prescribe the restriction 
\begin{equation}\label{eq:ghost_penalty_superdiagonal}
\eta_1\le 0\quad \text{and}\quad \eta_2\le 0.
\end{equation} Next, we consider the entry $\mathcal{B}_{GP}(k_1,k_1)$.This results in the condition
\[
\alpha h -\tau+\tau\eta_1+\tau\eta_2\overset{!}{\ge }0.
\]
Since $\alpha h -\tau$ is negative for $\alpha < \lambda$, we need to choose $\eta_1$ or $\eta_2$ to be positive. This is a contradiction to \eqref{eq:ghost_penalty_superdiagonal}. Therefore, it is \textit{not} possible to create a monotone scheme using this setup.

\subsubsection*{\newstablong}
We now consider the \newstab, which we introduced in \cite{May2019}.
The resulting scheme is of the same form as \eqref{eq:gp_formulation}, but instead of adding $J_h^\text{GP}$ we use the term 
\begin{equation}\label{eq:newstab}
J_h^{\newstababbrv}(u_h,w_h):=\beta\eta\jump{u_h}_{k-\frac{1}{2}}\jump{w_h}_{k_\text{cut}}.
\end{equation}
One big difference between \eqref{eq:gpstab} and \eqref{eq:newstab} is that the locations of the jump terms were moved. As a result, the position of the stabilization terms in the matrix $\mathcal{B}$ changed:  
\begin{equation*}
  \mathcal{B}_{\newstababbrv}=
    \left(
    \begin{smallmatrix}
      h-\tau & 0 & \cdots&  & & \cdots & 0 & \tau\\
      \tau & h-\tau & 0 & & & & & 0\\
      0 & \ddots & \ddots & & & & &\vdots\\
      \vdots & & \tau\left(1-\eta\right)& \diagentry{\alpha h-\tau\left(1-\eta\right)} & 0 & & &\\
      & & \tau\eta& \tau\left(1-\eta\right)& \diagentry{(1-\alpha)h-\tau} & & &\vdots \\
      \vdots & & & & \ddots& \ddots& & 0\\
      &&&&&&&\\
      0 &\cdots& & &\cdots&0&\tau& h-\tau\\
    \end{smallmatrix}\right).
\end{equation*}
In \cite{May2019}, we examined the monotonicity conditions of the stabilized scheme for the theta-scheme in time and a fixed value of $\eta$. Here, we will focus on
using explicit Euler in time and vary $\eta$ instead. Requiring that all entries
become non-negative results in the following three inequalities:
\begin{alignat*}{2}
\textbf{I}& \qquad  \alpha h -\tau(1-\eta) &\quad\ge 0,\\
\textbf{II}& \qquad  \tau \eta &\quad\ge 0,\\
\textbf{III}& \qquad  \tau(1-\eta) &\quad\ge 0.
\end{alignat*}
 Short calculations show that this implies the following restrictions on $\eta$
%\[
%0\overset{\textbf{II}}{\le} 1-\frac{\alpha}{\lambda}\overset{\textbf{I}}{\le} \eta
%\overset{\textbf{III}}{\le} 1, \quad \text{i.e., we need to choose} \quad \eta \in %\left[1-\frac{\alpha}{\lambda},1\right]
%\]
\[
\eta \overset{\textbf{II}}{\ge} 0, \quad
1-\frac{\alpha}{\lambda}\overset{\textbf{I}}{\le} \eta
\overset{\textbf{III}}{\le} 1, \quad \text{i.e., we need to choose} \quad \eta \in \left[1-\frac{\alpha}{\lambda},1\right]
\]
and should not stabilize for $\alpha > \lambda$. 
In other words, for $\alpha \ll \lambda < \frac{1}{2}$, the resulting scheme using
explicit Euler in time is monotone for $\eta \in \left[1-\frac{\alpha}{\lambda},1\right]$, despite the CFL condition on the cut cell $k_1$ being
violated. 
Next, we will discuss how to best choose $\eta$ within the prescribed range.
\section{Choice of $\eta$ in \newstab}\label{sec:4}
We denote the discrete solution on cell $j$ at time $t^n$ by $u^n_j$. Resolving the system
$\mathcal{M} u^{n+1} = \mathcal{B}_{\newstababbrv} u^n$ for the update on the two cut cells under the condition $\alpha < \lambda < \frac{1}{2}$, 
we get
\begin{align*}
u^{n+1}_{k_1} &= u^n_{k_1}-\frac{\lambda}{\alpha}\left(1-\eta\right)\left(u^n_{k_1}-u^n_{k-1}\right),\\
u^{n+1}_{k_2} &= u^n_{k_2}-\frac{\lambda}{1-\alpha}\left(u^n_{k_2}-u^n_{k_1}\right)-\frac{\lambda}{1-\alpha}\eta\left(u^n_{k_1}-u^n_{k-1}\right).
\end{align*}
For monotonicity, we need to choose $\eta\in \left[1-\frac{\alpha}{\lambda},1\right]$.
We will now examine the two extreme choices, $\eta=1-\frac{\alpha}{\lambda}$ and $\eta =1$, in more detail.

For $\eta=1-\frac{\alpha}{\lambda}$, the two update formulae reduce to 
\[
u^{n+1}_{k_1}=u^n_{k-1} \quad \text{and}\quad  u^{n+1}_{k_2} =\left(1-\frac{\lambda}{1-\alpha}\right)u^n_{k_2}+\frac{\alpha}{1-\alpha}u^n_{k_1}+\frac{\lambda-\alpha}{1-\alpha}u^n_{k-1}.
\]
We observe, comparing with figure \ref{fig:MP}, that the new update formulae now use the correct domains of dependence. In particular, $u^{n+1}_{k_1}$ now coincides with $u_{k-1}^n$ and $u^{n+1}_{k_2}$ now includes information from $u^n_{k-1}$, which is the neighbor of its inflow neighbor. Actually, the resulting updates correspond to exactly advecting a piecewise constant solution at time $t^n$ to time $t^{n+1}$ and to then averaging.
Therefore, thanks to the stabilization, we have implicitly restored the correct domains of dependence. In that sense, the new stabilization has a certain similarity to the
$h$-box method \cite{hbox_2003}.

For the choice $\eta=1$ the update formulae have the following form:
\[
u^{n+1}_{k_1} = u^n_{k_1} \quad \text{and} \quad u^{n+1}_{k_2} = u^n_{k_2} - \frac{\lambda}{1-\alpha}\left(u^n_{k_2}-u^n_{k-1}\right).
\]
We observe that in this case the smaller cut cell $k_1$ will not be updated. Instead, it just keeps its old value. In addition, the update of the solution on cell $k_2$ does not include information of its inflow neighbor $k_1$. Therefore choosing $\eta=1$ can be interpreted as skipping the small cut cell and let the information flow directly from its inflow neighbor into its outflow neighbor. 
\begin{remark}
In \cite{May2019}, we propose to use $\eta=1-\frac{\alpha}{2\lambda} \in \left[1-\frac{\alpha}{\lambda},1\right]$ as this produces better results for piecewise linear polynomials than $\eta=1-\frac{\alpha}{\lambda}$.
\end{remark}
\FloatBarrier
\section{Numerical results}
We will now compare the different choices of $\eta$ for the \newstab
numerically. We consider the grid described in figure \ref{fig:MP} and place cell $k$ such that $x_{k-\frac{1}{2}}=0.5$. 

\noindent
\begin{minipage}{\linewidth}
\noindent
We use discontinuous initial data
\begin{equation}\label{eq: disc init data}
    u_0(x) = \begin{cases} 1 & 0.1 \le x \le 0.5, \\
    0 & \text{otherwise},
    \end{cases}
\end{equation}
\end{minipage}
with the discontinuity being placed right in front of the small cut cell $k_1$.
We set $\beta = 1$, $\alpha = 0.001$, $\lambda = 0.4$, and $h=0.1$, and use $V_h^0(I)$ as well as periodic boundary conditions.
We test four different values for $\eta$: the extreme cases $\eta = 1$ and 
$\eta = 1-\frac{\alpha}{\lambda}$ as well as $\eta=1-\frac{\alpha}{2\lambda}$ and a value,
$\eta=1-\frac{2\alpha}{\lambda}$, that violates the monotonicity considerations.

\begin{figure}
\centering
\hspace*{\fill}%
%\subfloat[Test]
{\includegraphics[width=.46\linewidth]{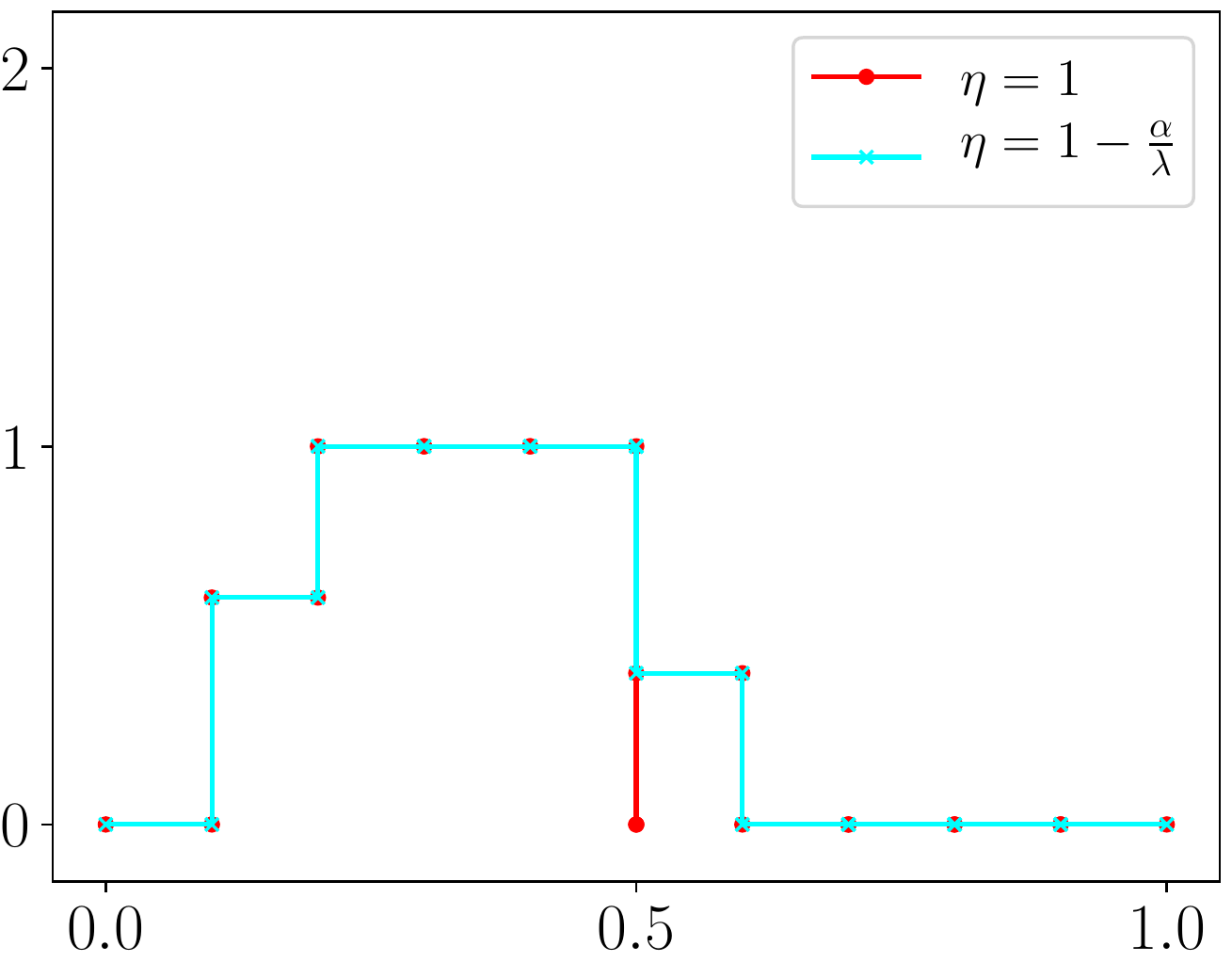}}%
\hfill
%\subfloat[Test]
{\includegraphics[width=.46\linewidth]{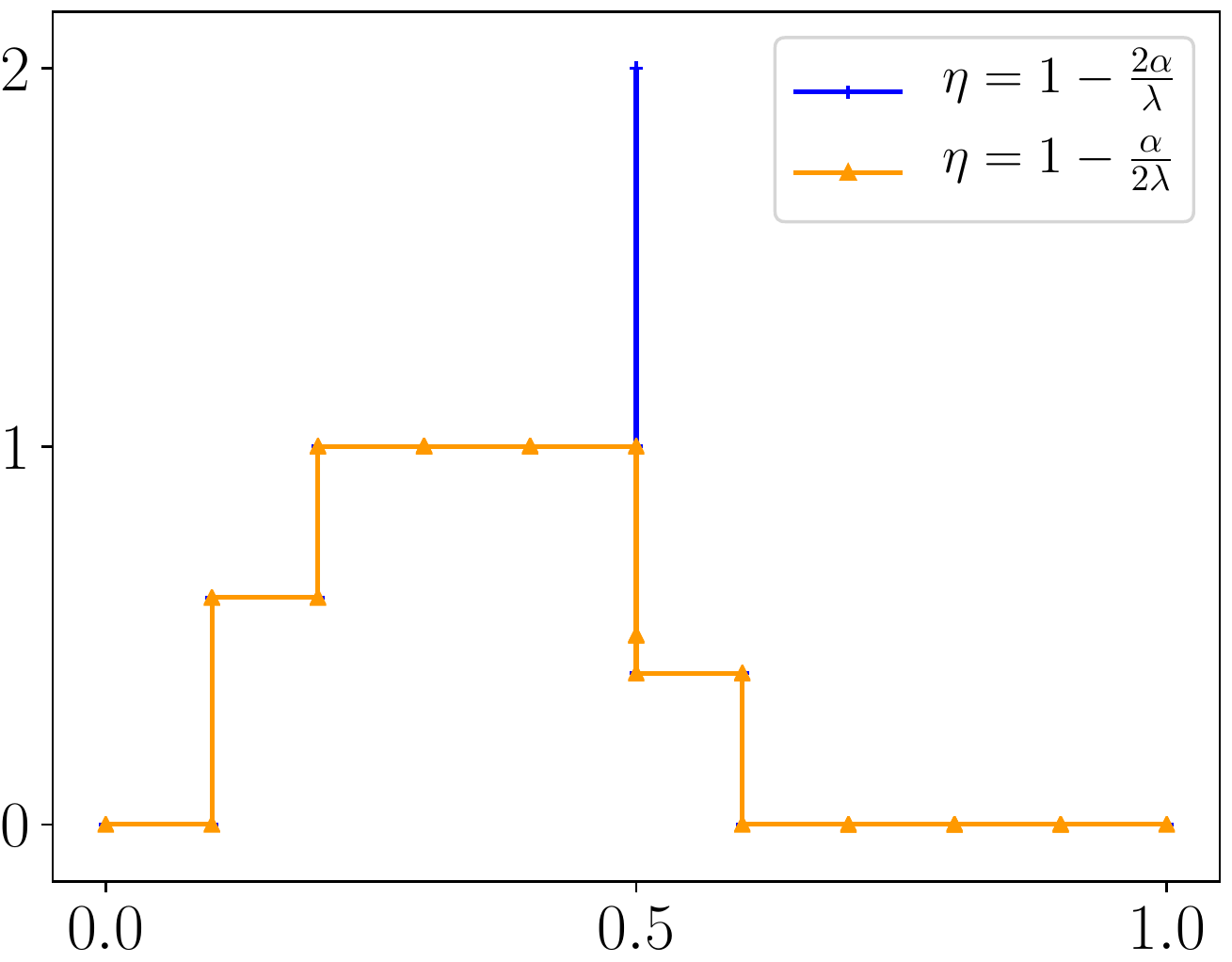}}%
\hfill
\caption{Results after one time step for the \newstab for different values of $\eta$. }\label{fig:disc}
\end{figure}
In figure \ref{fig:disc} we show the different solutions after one time step. 
For $\eta=1$ we observe that the solution on cell $k_1$ has not been updated, while the updates on the other cells are correct. Obviously, cell $k_1$ has simply been skipped.
The solution for $\eta=1-\frac{\alpha}{\lambda}$ corresponds to exactly advecting the initial data and to then apply averaging. 
If we choose $\eta=1-\frac{\alpha}{2\lambda}$, we observe that $u^1_{k_1}$ lies between $u^1_{k-1}$ and $u^1_{k_2}$.
Finally, for $\eta=1-\frac{2\alpha}{\lambda}$, which is not included in the suggested interval,
we observe a strong overshoot on the small cut cell. This cannot happen for a monotone scheme.
Therefore, the numerical results confirm our theoretical considerations above.
\FloatBarrier

\end{document}